\documentclass{article}
\usepackage{latexsym}
\usepackage{amsmath}
\usepackage{amsfonts}
\newtheorem{theorem}{Theorem}
\newtheorem{definition}{Definition}
\newtheorem{notation}{Notation}
\newtheorem{corollary}{Corollary}
\newtheorem{lemma}{Lemma}
\newtheorem{remark}{Remark}
\newtheorem{example}{Example}
\newcommand{\field}[1]{\mathbb{#1}}
\newcommand{\C}{\field{C}}
\newcommand{\R}{\field{R}}
\newcommand{\N}{\field{N}}
\newcommand{\Z}{\field{Z}}
\numberwithin{equation}{section}
\numberwithin{lemma}{section}
\numberwithin{theorem}{section}
\numberwithin{corollary}{section}
\numberwithin{remark}{section}
\numberwithin{definition}{section}

\begin{document}

\title{Distribution of the $k$-Multiple Point Range in the Closed Simple Random Walk I}

\author{Daniel H\"of\\
WorldLight Technologies \\
\" Ozer \.I\c shan\i {} K.:2 No.: 227\\
Ertu\u grul Bey Sok. 5 \\
\.Istiklal Mah. \\
TR-26010 Eski\c sehir \\
Republic of Turkey}
\maketitle

\begin{abstract}
The distribution of the number of points of the closed simple random walk, visited a given number of times (the $k$-multiple point range) is ana-lysed by a graph based approach. A general expression for the moments is derived. In this paper the joint generating function for dimension one is completely calculated and analysed for large lengths.
\end{abstract}

 \hfil \\

 \hfil \\
\section{Introduction}\label{S:intro}
In this paper we discuss the joint distribution of the variable $N_{2k}(w)$,
the number 
of points of a closed simple random walk $w$ on $\Z^{d}$ visited by
$w$ exactly $k$ times (the $k$-multiple point range).

In the vast literature about the $k$-multiple point range of
random walks \cite {Mard,Tek} the first moment of the 
distribution has 
been calculated and evaluated in the limit of large length of the walks. From this
Pitt \cite{Pitt} has obtained a law of large numbers for 
$N_{2k}(w)$ in the
case of transient random walks. Hamana in a series of papers \cite {Ham1, Ham2}
has 
also analysed the second moment of the distribution in the limit
of large lengths and obtained a central limit theorem for transient
random walks. For the recurrent simple random walk in two dimensions
Flatto \cite {Fla} has calculated the asymptotic first moment of the distribution and 
proven a law of large numbers. Hamana has analysed the asymptotic second
moment of the distribution for twodimensional recurrent random walks \cite {Ham3} and proven theorems about the limit distribution in this case. 
The distribution of the range, i.e. the number of different sites visited by the random walk has 
also been analysed in the literature \cite {Jain} up to the second moment. 

In this paper in section \ref{S:dis} we rederive the results about the first moment
of the distribution by a graph theoretical
method in the
context of the much more general theorem \ref{T:mom} about the algebraic structure of all
moments of the distribution.

In section \ref{S:one} we calculate the joint distribution of the $k$-multiple point range
for $d = 1$ from
these formulas in theorem \ref{T:one}. It turns out to be a generalized
geometrical distribution.
As simple applications we give the distribution of the
number of singlepoints, $N_2(w)$, in equation \ref{E:dsin} 
and the characteristic function of the distribution of doublepoints, 
$N_4(w)$, in equation \ref{E:ddou} and the distribution of the range
of a onedimensional simple
closed random walk in equation \ref{E:dvol}.

In section \ref{S:asy} we derive a method to calculate
the asymptotic (large length)
behaviour of the
joint distribution for the onedimensional closed simple walk. As a
typical application we obtain the asymptotic distribution of $N_2(w)$ and 
$N_4(w)$ in equation \ref{E:dou} and the asymptotic
joint second moment of the
$k$-multiple point range in equation \ref{E:2mom}. The asymptotic moments of the
distribution of the
range normed by it's asymptotic mean value are also given in equation 
\ref{E:avol} and turn out
to be the integer values of the $\xi$ - function of Riemann.
In section \ref{S:out} the results are discussed and an outlook on further
generalisations is given.
\section{ Distribution of moments}\label{S:dis}
\subsection{ Notations}
Throughout this article we work on the $d$-dimensional \emph{hypercubic lattice}
$L = \Z^{d}$. We denote the unit vector pointing in the $i$'th direction
by $e_i$, it's negative by $e_{-i} := -e_i$. \\

Because of our geometrical approach
we do not use the standard notation for a simple random walk. Instead we 
choose a \emph{geometrical definition of individual walks,} the corresponding
probability distributions are then obtained by summing variables over all walks e.g. of 
a given length 
and dividing by the number of walks of this length. Of course this is an equivalent 
formulation for the simple random walk. \\

So we denote a walk  by a pair $w = (p,s)$ with the starting point $p \in L$,
$p = (p^{(1)},\ldots,p^{(d)}) $
and a
sequence $s := (s_1,\ldots,s_n)$  of steps
$s_i \in \Z-0; -d \leq s_i \leq d$. The \emph{$i$'th point} of this walk
is then $p_i := p + \sum_{j=1}^{i} e_{s_{j}}; \; p_{0} := p$;
$i$ is called the \emph{index} of $s_i$ in $s$, $p_0$ is called the \emph{starting,}
$p_n$ the \emph{ending point.} The number $n$ is called the
\emph{length} of the walk $w$, i.e. $n = length(w)$. \\
A walk is called \emph{closed}
iff $p = p_n$. \\

Throughout this paper we use the term \emph{graph} in the meaning of \emph{ordered
digraph with loops.} 
Therefore we denote a graph $G$ by a triple $G=(V,E,J) $ consisting of a finite set 
of vertices $V$, a finite set of edges $E$ and a map $J$  giving the 
starting point $J(e)_1 \in V$ and ending point $J(e)_2 \in V $ of each edge $e$. The
edge e is then said to be an \emph{ingoing edge} of $J(e)_2$ and an \emph{outgoing edge}
of $J(e)_1$. \\

For $x,y$ elements of a discrete set we denote the \emph{Kronecker symbol} $\delta_{x,y}$ to be
$1$ for $x = y$ and $0$ otherwise. \\
The symbol $ \N $ is supposed to include $0$. \\ 
The symbol $\sharp S$
denotes the number of objects in a (finite) set $S$.\\
The notation $S_r$ denotes the
group of all permutations on r objects. \\

For $x,y \in \Z;\; y \neq 0$ the notation 
$x\%y$ denotes the nonnegative \emph{remainder} of the division of $x$ by $y$. \\ 

The symbol $z\diamond k$ is defined as follows:
If $g(z) := \sum_{n=0}^{\infty} a_n \cdot z^n$ is a formal series in $z$ with complex
coefficients then for any $k \in \N$ we denote 
$(g(z))_{z \diamond  k} := \sum_{n=0}^k a_n \cdot z^n$, i.e. the projection of the 
series onto the ring $R_k$ of polynomials of degree $k$.

\subsection{Archetypes}
In this subsection we establish a natural connection between 
graphs and walks.
\begin{definition}
Let $q \in L$  and $w = (p,s)$  a walk of length $n$ The multiplicity
$mu(q,w)$ of $q$ with respect to $w$ is defined as  $mu(q,w)
= \delta_{q,p} + \delta_{q,p_n} + 2\cdot\sum_{i=1}^{n-1} \delta_{q,p_i} $
, i.e. the number of incoming steps plus the 
number of outgoing steps.
\end{definition}
\begin{definition}
For a closed walk $w$ and an integer $k \in \N - {0}$ we define
the number $N_{2k}(w)$, the $k$-multiple point range as
$N_{2k}(w) := \sharp \{q \in L \mid mu(q,w) = 2k \} $ We also define the range
$ran(w) := \sum_{k=1}^{\infty } N_{2k}(w) $ 
\end{definition}
\begin{definition}
Let  $G=(V,E,J)$ be a graph and $ v \in V $  a vertex. Then 
the degree $deg(v,G)$ is defined as the number of ingoing edges 
plus the number of outgoing edges. 
\end{definition}
\begin{definition}
Let $ w=(p,s) $ be a closed walk of length $ 2n $. Let $ u = (q,t) $
be a (not necessary closed) walk of length $ m \leq 2n $ such that there
exists an index $ 1 \leq ind \leq 2n $ such that 
$ i = 1,\ldots,m \Rightarrow  t_{i} = s_{1+(i + ind -1)\%(2n) }$ and 
$ q = p_{ind} $. Then the pair $ (ind;u) $ is called a fixed subwalk of $ w $
\end{definition}
\begin{definition}
Let $(ind_1;u_1)$ and $(ind_2;u_2)$ be two fixed subwalks of a closed walk $w = (p,s)$
with lengths $m_1$ and $m_2$ and $2n$ and sequences $s^{(1)}$, $s^{(2)}$ and $s$ respectively. 
If $1 + (ind_1 + m_1 -1)\%(2n) = ind_2$ then the two fixed subwalks can be concatenated to 
the fixed subwalk $(ind_1;u_3)$ with $u_3 = (p_{ind_1},s^{(3)})$ with 
\begin{equation}
s^{(3)}_i := 
\begin{cases}
s^{(1)}_i \Longleftarrow i=1,\ldots,m_1 \\
s^{(2)}_{i-m_1} \Longleftarrow i=m_1+1,\ldots,m_1+m_2
\end{cases}
\end{equation}
\end{definition}
\begin{definition}
Let $(ind_1;u_1)$ and $(ind_2;u_2)$ be two fixed subwalks of a closed walk $w = (p,s)$
with lengths $m_1$ and $m_2$ and $2n$ and sequences $s^{(1)}$, $s^{(2)}$ and $s$ respectively. The fixed subwalks are said to overlap if there are numbers 
$1 \leq i_1 \leq  m_1, 1 \leq i_2 \leq m_2$ such that 
$(ind_1 + i_1)\%(2n) = (ind_2 + i_2)\%(2n)$
\end{definition}
\begin{definition}
Let $w=(x,s)$ be a closed walk of length $2n$. A refinement  of  $w$ is 
a collection of fixed subwalks of $w$ which do not overlap and can be concatenated to yield a fixed subwalk of length $2n$, i.e. loosely speaking whose concatenation
yields $w$.
\end{definition}
\begin{definition}\label{D:Arch}
Let $\Xi := (G_1,\ldots,G_k) $  be a family of $k$ Eulerian graphs \\
$ G_j = (V_j,E_j,J_j) $
with mutually 
disjoint edge sets. Let $\Sigma := (w_1,\ldots,w_k)$ be a family of closed walks and   
$P_j$ the set of points of $w_j$  with nonzero multiplicity and $S_j$ the set 
of fixed subwalks of $ w_j $ . A $2k$ tupel of maps
$ f = (f_{V_1},\ldots,f_{V_k};f_{E_1},\ldots,f_{E_k})  $  
is called \textbf{ archetype }
if and only if 
\begin{description}
\item[Structure] $ f_{E_i}:E_i \mapsto S_i  $ 
is injective and maps the different edges 
of $ E_i $  onto a refinement of $ w_i $ . $ f_{V_i}:V_i \mapsto P_i $  is injective 
and for $ v \in V_i \cap V_j \Rightarrow f_{V_i}(v) = f_{V_j}(v) $ .
\item[Orientation] $f_{V_i}$ maps the starting and ending vertices of edges as 
given by $J_i$ on 
the starting and ending points of the corresponding fixed subwalks.
\end{description}
\end{definition}
\begin{lemma}
With the notations of the above definition the following inequality is true  
\begin{eqnarray}
deg(v,G_j) \leq mu(f_{V_j}(v),w_j) \nonumber 
\end{eqnarray}
\end{lemma}
\emph{Proof.}
As $ G_j $  is Eulerian $ deg(v,G_j)/2 $  is the 
number of incoming edges and 
also the number of outgoing edges of $ v $.  $ f_{E_j} $ is injective, 
so any incoming and outgoing edge is mapped on a seperate subwalk (a loop being
mapped on one subwalk) which
starts or ends in $f_{V_j}(v)$. 
\subsection{Enumeration}
In this subsection we use the concept of an archetype to enumerate 
the moments of the distribution of points of a given multiplicity 
of a walk. 
\begin{definition}
We define the combinatorial symbol $P$ by the 
equation 
\begin{equation}
\prod_{k \in \N} (1 + \lambda_k)^{n_k} = 1 + \sum_{l=1}^{\infty }
\left(
\sum_{0 \leq k_1 \leq \ldots \leq k_l < \infty } P(n_{k_1},\ldots,n_{k_l}) 
\cdot \prod_{j=1}^l \lambda_{k_j}
\right)
\end{equation}
where $n_k \in \N$ and the $\lambda_k$ are variables of which only
finitely many are different from $0$.
\end{definition}
\begin{theorem} \label{T:fwa}
Let $w$ be a closed walk with $N_{2k}(w)$  points of multiplicity 
$2k$. Let $\Theta =\Theta (m_1,\ldots,m_l)$
be the set of equivalence classes under isomorphy of 
all Eulerian graphs with exactly $l$  vertices with degrees 
$2m_1,\ldots,2m_l$ . 
Let $ \sharp Arch(G,w) $ denote the number of archetypes from a Eulerian 
graph $ G $ into a walk $w$ . Let $\sharp Aut(G)$ denote the number of 
automorphisms of $ G $. Then
\begin{equation}\label{E:en1}
\sum_{[G] \in \Theta (m_1,\ldots,m_l)} \frac {\sharp Arch(G,w)} {\sharp Aut(G)} =
\sum_{k_i \geq m_i} P(N_{2k_1}(w),\ldots,N_{2k_l}(w)) \cdot 
\prod_{j=1}^{l} \binom {k_j} {m_j}
\end{equation}
where [G] denotes the class. 
\end{theorem}

\emph{Proof.}
Let $ (p^{[1]},\ldots,p^{[l]}) $
  be $ l $  different points of $ w = (p,s) $ with multiplicities 
$ 2k_1,\ldots,2k_l $
. This means that for each j separately there exist $ k_j $ different
 indices $ ind(j,\alpha^{(j)}) \in 
\{ 1,\ldots,length(w) \} $ with $ \alpha^{(j)} = 1,\ldots,k_j $
 such that $ p^{[j]} = p_{ind(j,\alpha^{(j)})} \quad 
\forall \alpha^{(j)} = 1,\ldots,k_j$. 
For each 
$ j $ there are $ \binom {k_j} {m_j} $
 ways to choose $m_j$  different indices $ ind(j,\alpha^{(j)}_q) $ with 
$ q = 1,\ldots,m_j $ . For such a choice we order 
the chosen indices of $ s $ in the form 
\begin{equation}\label{E:ind}
1 \leq ind^{(1)} < \ldots < ind^{(m)} \leq  length(w) 
\end{equation}
with 
$ m = \sum_{j=1}^{l} m_j $. For each such choice of points 
and indices  we now construct a different archetype. We define 
the graph $G=(V,E,J)$ by letting $V = \{ p^{[1]},\ldots,p^{[l]}\}$, $ E $ 
being 
the set of  fixed subwalks of $w$ which correspond to sequences 
between adjacent indices in our ordered set in $s$. $J$ then maps 
these subwalks on the corresponding starting and ending point. 
With this definition $G$ is obviously a graph and $[G] \in \Theta $, 
and we can define 
an archetype $f$ from $G$ to $w$ by the obvious projections. 
Now if $f$ is any archetype from a graph $G=(V,E,J)$ then $V$ can 
be identified with $f_V(V)$ and $E$ with $f_E(E)$ because of the injectivity 
of the respective maps. Now $f_E(E)$ is a refinement of $w $ and by using the
orientation property of $f_V(V)$ we get an ordered index system in the form 
of equation \ref{E:ind} with the corresponding starting and ending points according
to $J$. But from this it is obvious that $G$ is isomorphic to a graph in the
class we constructed previously. 

\begin{theorem}\label{T:fGr}
Let $G=(V,E,J)$ be a Eulerian graph with vertices $v_1,\ldots,v_l$
with degree $m_1,\ldots,m_l$. Let $Y :=(y_1,\ldots,y_l) \in L^l$ 
be a vector of $l$ different points in $L = \Z^d$
Let 
\begin{equation}
\begin{array}{rccl}
\pi:& V & \longmapsto & \{y_1,\ldots,y_l\} \\
    & v_i & \longmapsto & y_i 
\end{array}
\end{equation}
be the natural projection.\\ Let $g_{2n}(G,\pi,Y)$  denote the number of different archetypes 
$(f_V;f_E)$ of $G$ into closed walks of length $2n$ with $f_V = \pi$ and their generating function
\begin{equation}
N(z,G,\pi,Y) := \sum_{n=0}^{\infty } g_{2n}(G,\pi,Y)\cdot z^{2n}
\end{equation}
then the following equation 
\begin{equation} \label{E:fGr}
N(z,G,\pi,Y) = N_E(G)\cdot z \cdot \frac{\partial }{\partial z} 
\left( \prod_{e \in E} h(\pi(J(e)_1) - \pi(J(e)_2),d,z) \right)
\end{equation}
holds, where $N_E(G)$ is the number of Euler trails of $G$; for $p \in L$, 
the symbol $h$ stands for 
\begin{equation} \label{E:hxdz}
h(p,d,z) := -\delta_{p,0} + \int_0^{\infty }dy \cdot e^{-y} \cdot
\prod_{\alpha = 1}^d I_{\mid p^{(\alpha)} \mid}(2zy)
\end{equation}
with the modified Bessel function $I_{\mu } $  and 
$\mid z \mid < \frac{1} {2d} $
\end{theorem} 

\emph{Proof.}
If $(f_V;f_E)$ is an archetype from $G$ into a walk $w$ with 
$f_V = \pi$ then $f_E(E)$ is a 
refinement of $w$. We can order the fixed subwalks in this refinement 
in the form of equation \ref{E:ind} according to the value of their indices
in $w$. By the injectivity of $f_E$  
we can transfer this order of the subwalks to a unique order of the edges of $G$. 
From the orientation 
property of  $f_V$  it is clear that this order induces a unique 
Euler trail on $G$. 

Now let an Euler trail on $G$ be given by an ordered sequence 
of the edges $e^{(1)},\ldots,e^{(q)}$ . Then we know that any 
archetype $(f_V;f_E)$  in the form of the 
theorem belongs to a walk $w$ which visits the points $y_i$ in the 
sequence induced by the Euler trail and this sequence determines 
$(f_V;f_E)$ uniquely once the subwalks corresponding to the edges
and the position of the starting point in the collection of subwalks
of $w$ are given. 
Now as $h(p,z,d)$ is the generating function 
for the number of walks between the point $0$ and $p$ in $L$,  by 
the Markov feature the product   
\begin{displaymath}
z \cdot \frac{\partial }{\partial z} 
\left( 
\prod_{1 \leq \mu \leq q } h(\pi(J(e^{(\mu}))_1) - \pi(J(e^{(\mu)})_2),d,z) 
\right)
\end{displaymath}
is the generating function 
for the number of such collections and choices of the starting point.
But this factor is the same for all Euler trails of $G$ because the 
product in the brackets in fact is a product over all edges in $E$ which prooves 
the theorem. 
\begin{definition}
Let $r$ (not necessarily different) numbers $k_i \in \N - \{0\}; 
i = 1,\ldots,r $ be given. 
On the $r$-permutations $\sigma \in S_r$ we define an equivalence relation $\sim $
by 
\begin{equation}
\sigma \sim \hat \sigma \Longleftrightarrow 
k_{\sigma (i)} = k_{\hat \sigma (i)} \quad \forall i=1,\ldots,r
\end{equation}
Then we define $D(k_1,\ldots,k_r)$ to be the set of equivalence classes of 
permutations of $r$ objects under the above equivalence relation.
\end{definition}

\begin{definition}\label {D:hr}
Let there be  $r$ natural numbers $h_i \in \N - \{0\};i=1,\ldots,r$ and
$r \times r$ Matrices  $(F_{i,j}) \in GL(r,\N)$ which fulfil the following
conditions:
\begin{equation}\label{E:bal}
\begin{array}{c}
F_{i,i} := 0\\
\sum_{i=1}^r F_{i,j} = \sum_{i=1}^r F_{j,i} = h_i \\
\forall i = 1,\ldots,r
\end{array}
\end{equation}
\begin{equation}\label{E:cof}
cof(A-F) \neq 0
\end{equation}
where $A$ is the diagonal $r \times r$ Matrix $A = diag(h_1,\ldots,h_r)$ and cof is
the cofactor (A matrix fulfiling equation \ref {E:bal} has cofactors which are the same  no matter
which line or which column is erased) 
We define the set of Matrices $(F_{i,j})$ which fulfil only
equation \ref {E:bal} with $H_r(h_1,\ldots,h_r)$, and with $\tilde H_r(h_1,\ldots,h_r)$
the set of those which fulfil 
equations \ref {E:bal} and \ref {E:cof}. 
\end{definition}

\begin{theorem} \label{T:mom}
Let  $2k_1,\ldots,2k_r$ be $r$ (not necessarily different) multiplicities. 
Let $J_r := \{ Y:=(Y_0,\ldots,Y_{r-1}) \in L^r \mid Y_0 = 0 \wedge i \neq j \Longrightarrow Y_i \neq Y_j \}$.
Let $W_{2N}$ denote the set of closed random walks which start and end in $0$
and with $length(w) \leq 2N$. Then the moments of the distribution of $N_{2k}(w)$
are given by

\begin{multline}\label{E:mome}
\sum_{w \in W_{2N}} P(N_{2k_1}(w),\ldots,N_{2k_r}(w)) \cdot z^{length(w)} =  \\
\frac {1} {r!} \cdot z \cdot \frac 
{\partial } {\partial z}
\Bigl[ 
\sum_{\sigma \in D(k_1,\ldots,k_r)} 
( 
\sum_{Y \in J_r}
( 
\sum_{h_i = 1}^{\infty }
( 
\sum_{F \in H_r(h_1,\ldots,h_r)} \\ 
cof(A-F) 
\left( \prod_{\substack{1 \leq a \leq r \\ 1 \leq b \leq r}} \frac {h(Y_a - Y_b,d,z)^{F_{a,b}}} 
{(1+ h(0,d,z))^{F_{a,b}}\cdot F_{a,b}\,!}
\right) \cdot  \\
\prod_{j=1}^r (-1)^{h_j} \cdot h_j\,!\cdot K(h_j,k_{\sigma (j)},h(0,d,z))
)
)
)
\Bigr]_{z \diamond 2N}
\end{multline}
($z \diamond 2N$ see Notations) where the vertex factors $K$ are given by the equation
\begin{multline} \label {E:vert}
K(q,k,\omega ) :=  \\ 
\begin{cases}
(-1)^{k}\cdot \frac{1}{q} \cdot \left( \frac {1} {1 + \omega } \right)^k
\cdot \sum_{\nu = max(0,k-q)}^{k-1} \binom {k-1} {\nu} \cdot 
\binom {q} {k - \nu} \cdot (- \omega)^{\nu}& \Longleftarrow q > 0 \\
\frac {1} {k} \cdot \left( \frac {\omega }{1 + \omega} \right)^k
& \Longleftarrow q = 0 
\end{cases}
\end{multline}
\end{theorem}

\emph{Proof.}
From equation \ref{E:hxdz} it is immediately clear that the Taylor expansion of 
$h(x,d,z)$ in $z$ around the origin, for $x \neq 0$ 
starts with the power $z^{(\sum_{\alpha = 1}^d \mid x^{(\alpha)} \mid)}$.
So for $ \mid Y_a - Y_b \mid > 2N$ the factor 
$(h(Y_a - Y_b,d,z))_{z \diamond 2N} \equiv 0$ 
or for 
$\sum_{\substack{a=1 \\ b=1}}^r F_{a,b} > 2N $  the product

\begin{displaymath}
\left( \prod_{\substack{a=1 \\ b=1}}^r \frac {h(Y_a - Y_b,d,z)^{F_{a,b}}} 
{(1+ h(0,d,z))^{F_{a,b}}\cdot F_{a,b}\,!}
\right)_{z \diamond 2N} \equiv 0 
\end{displaymath}
as $a \neq b \wedge Y \in J_r \Longrightarrow 
\mid Y_a - Y_b \mid > 1$. 
So we do not have to deal with convergence problems
in equation \ref{E:mome}; because of the projection $ z \diamond 2N$, 
all sums are in fact finite.

For the proof of the identity we combine Theorems \ref{T:fwa} and 
\ref{T:fGr} in the following way:
We take equation \ref{E:fGr}, apply the operation $()_{z \diamond 2N}$ 
and then sum it over all possible vectors $Y \in J_r$. We insert this into 
equation \ref{E:en1}. Using translation invariance we find:
\begin{multline}\label{E:bir}
\sum_{\substack{Y \in J_r \\ [G] \in \Theta (m_1,\ldots,m_r)}}
\frac {N_E(G)} {\sharp Aut(G)}  
\cdot z \cdot \frac {\partial } {\partial z}
\left(
\prod_{e \in E} h(\pi_Y(J(e)_1) - \pi_Y(J(e)_2),d,z)
\right)_{z\diamondsuit 2N} \\
  = \sum_{w \in W_{2N}} \sum_{k_i \geq m_i} 
\left(
P(N_{2k_1},\ldots,N_{2k_r}) \prod_{i=1}^r \binom {k_i} {m_i}
\right)
z^{length(w)}
\end{multline}
where $\pi_Y$ denotes the projection of theorem \ref{T:fGr} specifically for 
the vector $Y$. 
An isomorphy class of a Eulerian graph $[G] \in \Theta (m_1,\ldots,m_r)$
fully determines an adjacency matrix $F \in \tilde H_r(h_1,\ldots,h_r)$ and 
loop numbers $l_i = m_i - h_i$ at the corresponding vertices. On the other 
hand the data $F \in \tilde H_r(h_1,\ldots,h_r)$ and $l_i$ fully
determine an isomorphy class of Eulerian graphs 
$[G] \in \Theta (h_1 + l_1,\ldots,h_r + l_r)$ up to a permutation of vertices.
So 
\begin{equation}
\sum_{[G] \in \Theta (m_1,\ldots,m_r)}  \longleftrightarrow 
\sum_{l_i = 0}^{m_i-1} \sum_{F \in \tilde H_r(m_1 - l_1,\ldots,m_r - l_r)} 
\frac {\sharp Aut(F)}{r!}
\end{equation}
where $Aut(F)$ is the group consisting of the permutations $\sigma \in S_r$ with 
the property $\sigma \in Aut(F) \Longleftrightarrow 
F_{\sigma(i),\sigma(j)} = F_{i,j} \quad \forall i,j$. \quad
In these new variables 
\begin{equation}
N_E(G) = cof(A-F) \cdot \prod_{i=1}^r (h_i - 1)!
\end{equation}
according to the generalized theorem of Tutte \cite {aar} and 
\begin{equation}
\sharp Aut(G) = \sharp Aut(F) \cdot \prod_{i=1}^r l_i! \prod_{k,j=1}^r F_{k,j}!
\end{equation}
To get rid of the combinatorial factors at the right hand side of 
equation \ref{E:bir} we note that for a sequence $(\gamma_k)_{k \in \N}$
of complex numbers of which only finitely many are different from zero the derived
sequence $g_m := \sum_{k \geq m} \gamma_k \cdot \binom {k}{m}$ again has only finitely
many nonzero entries and 
$ \gamma_m = \sum_{k \geq m} g_k \cdot \binom {k}{m} (-1)^{m+k} $
Putting all this together we arrive at equation \ref{E:mome}, noting that
the sum over $H_r$ instead of $\tilde H_r$ does not make a difference because of 
the factor $cof(A-F)$ in the sum and noting that the vertex factor 
\begin{equation}
K(q,k,\omega) := \frac {(-1)^q }{q!} (1 + \omega)^q \sum_{m \geq max(k,q)}
\binom {m} {k} \frac{ (m-1)! }{(m-q)!} (\omega)^{m-q} (-1)^{(m+k)}
\end{equation}
can be resummed in the form of equation \ref {E:vert}
\begin{corollary}
The first moment of the number of points hit 
exactly $2k$ times  by a closed  random walk $w$, starting in the 
point $0$ on the lattice $L= \Z^d$ is given by
\begin{equation}
\sum_{w} N_{2k}(w)\cdot z^{length(w)} = z \cdot \frac {\partial} {\partial z} 
\left( \frac {1} {k} \left( \frac {h(0,d,z)} {1 + h(0,d,z)} \right)^k \right)
\end{equation}
\begin {equation}
\sum_{w} ran(w) \cdot z^{length(w)} = z \cdot \frac{\partial} {\partial z}
\left( \log (1 + h(0,d,z))
\right) 
\end {equation}
where $h(0,d,z)$ is given in equation \ref{E:hxdz}. 
We use 
\begin{equation}
\begin{split}
h(0,1,z) &= \frac{1 - \sqrt{1 - 4z^2}} {\sqrt{1 - 4z^2}} \\
h(0,2,z) &= \frac {2}{\pi}K(16z^2) -1 \\
h(0,d,1/2d) &= G(d) := \int_0^{\infty} dy \cdot I_0(y/d)^d \cdot e^{-y} 
\end{split}
\end{equation}
where $K$ denotes the complete elliptic integral and the last equation is 
valid for $d > 2$.
We then have the following asymptotic (in terms of $length(w)$) behaviour 
of the first moment 
\begin {multline}\label{E:amo}
E_n(N_{2k}(w)) := \frac {\sum_{lenght(w)=2n} N_{2k}(w) } {\sum_{length(w) = 2n} 1} = \\
\begin {cases}
1 + O(1/n) & \Longleftarrow  d = 1 \\
\frac {2n  \cdot \pi^2} {log(n)^2}  + O(1/\log(n)) & \Longleftarrow  d = 2 \\
2n \cdot \frac {G(d)^{k-1} } {(1 + G(d))^{k+ 1}} + O(1,n^{2-d/2}) & \Longleftarrow  d > 2
\end {cases}
\end {multline}
reproducing the results from Mardudin et al. \cite{Mard} and Flatto \cite {Fla}
(note that $1/(1 + G(d))$ is the probability that the random walk never returns
to the starting point). 
For the range we find
\begin{equation}
E_n(ran(w)) = 
\begin{cases}
\sqrt{\pi n} \cdot (1 + O(1/n)) & \Longleftarrow d = 1 \\
\frac {2n \cdot \pi} {\log(n)} \cdot (1 + O(1/\log(n)) & \Longleftarrow d = 2 \\
2n \cdot \frac {1} {1 + G(d)} + O(1,n^{2 - d/2}) & \Longleftarrow d > 2 
\end {cases}
\end{equation}
\end{corollary}
In this paper we will concentrate on dimension one to keep the size of the paper limited. More results, especially higher moments in higher dimensions
will be published in later articles. For dimension one the following reformulation of equation \ref{E:mome} is useful. 
\begin {lemma}
Using the notations of theorem \ref{T:mom} the following formula is true
\begin {multline} \label{E:mdet}
\sum_{w \in W_{2N}} P(N_{2k_1}(w),\ldots,N_{2k_r}(w)) \cdot z^{length(w)} = 
 \\
\frac {1} {r!} \cdot z \cdot \frac {\partial } {\partial z} 
\Biggl[ 
\sum_{ \substack {Y \in J_r \\ \vartheta \in D(k_1,\ldots,k_r)}}
\left( \prod_{i=1}^r \hat K(x_i,k_{\vartheta (i)},h(0,d,z)) \right)
cof(\hat A - \hat F) \\
\left( \prod_{j=1}^r \frac {1} {(1-x_j)^2} \right)
\frac {1} {det(W^{-1} + X)} 
\Bigr\rvert_{\forall a: x_a = 0; \quad \forall b,c: X_{b,c}= U(Y)_{b,c}}
\Biggr]_{z \diamond 2N}
\end {multline}
with the notations
\begin{equation} \label {E:oper}
\begin{array} {rl}
\hat K(x,k,\omega ) :&= \left( \frac {1}{1+\omega } \right)^k \cdot
\sum_{\nu =0}^{k-1} \binom {k-1}{\nu } \frac {(\omega)^\nu } {(k-\nu )!} 
\left( \frac {\partial } {\partial x} \right)^{k- \nu -1} \\
\hat F_{i,j} :&= X_{i,j} \cdot \frac {\partial } {\partial X_{i,j}} \\
\hat A_{i,j} :&= \sum_{k=1}^r \hat F_{i,k}\\
W_{i,j} :&= \delta_{i,j} \cdot (1-x_i) \\
U(Y)_{i,j} :&= \frac {h(Y_i - Y_j,d,z)} {1 + h(0,d,z)} \cdot (1 - \delta_{i,j})
\end{array}
\end{equation}
and the cofactor $cof$ understood as the corresponding formal expansion in the differential operators.
\end {lemma}
\emph{Proof.}
For a complex $r \times r$ matrix $X$ with zero diagonal and with eigenvalues which have an absolute value
smaller than $1$, in the notations of \cite {QFT}
\begin{equation} \label{E:det}
\begin{split}
\frac {1} {det(1-X)} &= \frac {1} {(2\pi)^r} \int_{\R^{2r}} 
\prod_{i=1}^r d\varsigma_i\cdot d\bar \varsigma_i \cdot         
exp\left( - \sum_{j,k = 0} ^r \bar \varsigma_j (1 - X)_{j,k} \varsigma_k \right)
  \\
&= \sum_{h_i=0}^{\infty } \sum_{F \in H_r(h_1,\ldots,h_r)} 
\left( \prod_{\substack{a = 1 \\ b = 1}}^r  \frac {X_{a,b}^{F_{a,b}}} {F_{a,b}!}\right)
\left( \prod_{j=1}^r h_j! \right)
\end{split}
\end{equation}
holds where the second line follows from the expansion of the exponential in the
variables $X_{a,b}$ and subsequently integrating the Gaussian over the polynomials 
in the variables $\bar \varsigma_i, \varsigma_i$. Now the lemma follows from the 
right hand side of equation \ref{E:det}
by observing that the variables $X$ are conjugated to the variables $F$ and that 
\begin{equation}
\hat K(x,k,\omega)(x-1)^{q - 1} \Bigr\rvert_{x=0} = (-1)^{q+1}K(q,k,\omega)
\end{equation}
\section {The onedimensional case}\label{S:one}
\begin{lemma}\label{L:dif}
Let  f be a complex valued absolutely summable function on $J_r$  (as defined  in Theorem \ref{T:mom}) which is invariant under permutations of the components of 
$Y \in J_r$  and  only depends on the absolute value of the differences 
$ \mid Y_a - Y_b \mid $
of the components. For $ d = 1 $ we then have 
\begin {equation}
\frac {1} {r!} \sum_{Y \in J_r} f(Y) = \sum_{\substack {g_0 = 1 \\
\ldots \\
g_{r-2} = 1}}^{\infty } f(\tilde Y) \Bigr\rvert_{\tilde Y_0 := 0; 
\tilde Y_a := \sum_{j=1}^{a-1} g_j}
\end {equation}
\end{lemma}
\emph{Proof.}:  For $d = 1$ and $Y \in J_r$ there is exactly one permutation 
$\sigma  \in S_r $ of the indices  of $Y$  under which they are ordered according to their value in $\Z$.  Therefore we can write the components in the form   
$Y_{\sigma(a)} = Y_{\sigma(0)} + \sum_{j=1}^{a-1} g_j $
with a strictly positive $ g_j := Y_{\sigma(j+1)} - Y_{\sigma(j)} $.  As the function only depends on the differences of the components the lemma is true.
\begin {lemma}\label{L:perm}
For  $d = 1$ and any vector $Y \in J_r$ with $a < b \Longrightarrow Y_a < Y_b$  
and any matrix  $F \in H_r(h_1,\ldots,h_r)$ (see equation \ref {E:bal}) and 
any matrix $U(Y)_{i,j}$ in the form of equation \ref {E:oper} the following formula is true:
\begin {equation}
\prod_{\substack {i=1 \\ j = 1}}^r U(Y)_{i,j}^{F_{i,j}} =
\prod_{\substack {i=1 \\ j = 1}}^r \tilde U(Y)_{i,j}^{F_{i,j}} 
\end {equation}
with 
\begin{equation}
\tilde U(Y)_{i,j} := \begin{cases}
U(Y)_{i,j}^2 & \Longleftarrow  i \leq j  \\
1 & \Longleftarrow  i > j
\end{cases} 
\end {equation} 
\end {lemma}
\emph{Proof.}
For  $d = 1$  we have 
\begin{equation}\label{E:qkz}
\begin{split}
U_{i,j}(Y) & =  \prod_{k = min(i,j)}^{max(i-1,j-1)} q_k(z) \\
q_k(z) & = \left( \frac {2z} {1 + \sqrt{1-4z^2}} \right)^{g_{k -1}}
\end{split}
\end{equation}
Therefore 
\begin {equation}
\prod_{i,j=1}^r U_{i,j}(Y) = \prod_{k=1}^{r-1} 
q_k(z)^{ \sum_{\substack {l>k \\ m \leq k}} F_{m,l} + 
\sum_{ \substack {l \leq k \\ m > k}} F_{m,l}} 
\end  {equation} 
where the sums in the exponential refer to the upper and lower half of the matrix 
$U_{i,j}(Y)$ respectively.  Using equation \ref{E:bal} one easily sees that the two sums in the exponent on the right hand side are actually the same which proves the lemma. 
\begin {theorem}\label{T:one}
Let $(\Lambda_k)_{k \in \N-\{0\} }$ be a set of variables 
$\Lambda_k \in \C \quad \forall k$ and let us denote 
$t_k := e^{\Lambda_k} - 1$. Let 
$A(z) := \sqrt {1 - 4z^2}$ and $B(z) := \frac {2z} {1 + A(z)}$
The joint generating function for the distribution of the $k$-multiple point range in one 
dimension is given by
(generalized geometrical distribution)
\begin{multline}
\sum_{w \in W_{2N}} \left( exp \left( \sum_{k=1}^{\infty } \Lambda_k N_{2k}(w) 
\right) \right) \cdot 
z^{length(w)} = \\ 
z \frac {\partial } {\partial z} 
\Biggl[ 
T_0(z) + \sum_{k=1}^{\infty} t_k \cdot T_1(z)^{(k)} + 
\sum_{\substack {k_1 = 1 \\ k_2 = 1}}^{\infty} t_{k_1} \cdot t_{k_2} \cdot T_2(z)^{(k_1),(k_2)} + \\
\sum_{\substack {k_1 = 1 \\ k_2 = 1 \\ k_3 = 1}}^{\infty} t_{k_1} \cdot t_{k_2} \cdot t_{k_3} \cdot T_{\geq 3}(z)^{(k_1),(k_2),(k_3)} 
\Biggr]
\end {multline}
with the different terms $T$ given by 
\begin {equation}
\begin {split}
T_0(z) &:= \left( log \left( \frac {2} {1 + A(z)} \right) \right)_{z \diamond 2N} \\
T_1(z)^{(k)} &:= \left( \frac {1} {k} \left( \frac {1 - A(z)} {A(z)} \right)^k 
\right)_{z \diamond 2N}
\end {split}
\end {equation}
If we define 
\begin{equation}
G_{i,j}(z) := 
\left(
(-1)^{i+j}\cdot A(z)^{i+j} \sum_{f=max(i,j)}^{\infty} \frac {1} {f}
\binom {f} {i} \binom {f} {j} \frac {B(z)^{2f}} {1 - B(z)^{2f}}
\right)_{z \diamond 2N}
\end{equation}
Then 
\begin{multline}
T_2(z)^{(k_1),(k_2)} := \\ \left(
\sum_{\substack{0 \leq l_1 \leq k_1 -1 \\ 0 \leq l_2 \leq k_2 -1}}
\binom {k_1 - 1} {l_1} \binom {k_2 -1} {l_2} (1 - A(z))^{l_1 + l_2}
G_{k_1 - l_1,k_2 - l_2}(z) \right)_{z \diamond 2N}
\end{multline}
For the last term $T_{\geq 3}(z)^{(k_1),(k_2),(k_3)}$ we define index pairs $(\varrho,t); \quad 0 \leq \varrho,t < N$ for vectors 
$ \Phi(k_1,k_2,z), \Psi(k_3,z) \in \R^{N^2}$ whose components are given 
by
\begin {equation}\label {E:psid}
\begin {split}
\Psi_{(\varrho,t)}(k,z) &:=  \Biggr( \delta_{t,0} (-1)^{\varrho + 1} 
\sum_{l=0}^{k - 1} \binom {k-1} {l} (1 - A(z))^{l} G_{\varrho + 1,k - l}(z) 
\Biggl)_{z \diamond 2N} \\
\Phi_{(\varrho,t)}(k_1,k_2,z) &:=  \Biggr( \frac {(-1)^{\varrho + 1} \cdot (k_1 - 1)!} {\varrho ! (k_1 - 2 - \varrho - t)!} T_2(z)^{(k_1 - 1 - \varrho -t),(k_2)} \Biggl)_{z \diamond 2N}
\end {split}
\end {equation}
and matrices $Q(k) \in Gl(N^2,\R)$ with 
\begin {multline}\label {E:qdef}
Q_{(\varrho,t),(\tilde \varrho,\tilde t)}(k,z) := \Biggl( \frac {1} {\tilde t ! \cdot (\tilde t + 1)! \cdot \varrho !} \\  
\sum_{\eta =0}^{k-2-\varrho - t - \tilde t}
\frac {(k-1)! (1-A(z))^{k-2-\varrho - t - \tilde t - \eta}  (-1)^{\varrho + \eta}}
{\eta ! (k - 2 - \varrho - t - \tilde t - \eta)!} 
H_{\tilde t + 1,\tilde \varrho + \eta + 2 \tilde t + 2} (z) \Biggl)_{z \diamond 2N}
\end {multline}
with the notation
\begin {equation}
H_{i,j}(z) := \left( A(z)^j \sum_{f = j - i} \binom {f + i - 1} {j - 1} 
\frac {B(z)^{2f}} { 1 - B(z)^{2f}} \right)_{z \diamond 2N}
\end {equation}
With these notations
\begin {equation}
T_{\geq 3}(z)^{(k_1),(k_2),(k_3)} := 
\Biggl< \Psi(k_1,z)\Bigr\rvert(\mathbf{I} - \sum_{k = 1}^{\infty}t_k \cdot Q(k,z))^{-1} \star \Phi(k_2,k_3,z) \Biggr>_{z \diamond 2N}
\end {equation}
where $\mathbf{I}$, and $()^{-1}$ denote the unit 
element and the inversion in $Gl(N^2,\R)$, $ \star $ and $< . \mid . > $, 
the application of a matrix on a vector and the standard 
scalar product respectively in $\R^{(N^2)}$. Because of 
$z \diamond 2N$ the inversion
is actually a shorthand notation for the corresponding geometrical series. 
\end {theorem} 
\emph {Proof.}
We recall equation \ref{E:mome} and use lemma \ref{L:dif} for the sum over 
$Y \in J_r$  and then lemma \ref{L:perm} to replace $U(Y)$  by 
$\tilde U(Y)$. We transform the sum over the $r \times r$  matrices 
$F \in H_r$ into a complex Gaussian integral as in equation \ref{E:det} and can 
then also express the cofactor in terms of a determinant. As a result of these steps we get:
\begin {multline}
\sum_{w \in W_{2N}} P(N_{2k_1}(w),\ldots,N_{2k_r}(w)) \cdot z^{length(w)} = \\
z \frac {\partial}{\partial z} \Biggl(
\sum_{\substack {\vartheta  \in Da(k_1,\ldots,k_r \\ g_k \in \N-\{0\})}}
\prod_{i=1}^r \hat K(x_i,k_{\vartheta (i)},h(0,1,z))
\\
\frac {(-1)} {(2\pi)^r} \int_{\R^{2r}} \prod_{j=1}^r 
\frac {d\varsigma_j \cdot d\bar \varsigma_j}{(1 - x_j)^2} det(M(Y)) \cdot 
exp\left(
- \sum_{\substack {a=1 \\ b=1}}^r \bar \varsigma_a \left[ W^{-1}_{a,b} + 
\tilde U_{a,b}(Y) \right] \varsigma_b
\right) \Biggr)
\end {multline}
with the $(r-1)\times (r-1)$  matrix $M(Y)$
\begin {multline}
M_{a,b}(Y) := \begin {cases} 
\bar \varsigma_a \tilde U_{a,b}(Y) \varsigma_b &  \Longleftarrow  a \neq b \\
-\varsigma_a \sum_{c=1}^r \tilde U_{a,c}(Y) \varsigma_c & \Longleftarrow a = b
\end {cases} \\
(a=2,\ldots,r; \quad b=1,\ldots,r-1)
\end {multline}
and  the matrix $W$  defined as in equation \ref {E:oper} and the differences   
$g_k$ as in lemma \ref {L:dif}. 
By extracting the common row factors $\bar \varsigma_a$ out of the determinant, then adding an appropriately scaled first column to each of the other columns and permuting the columns, $M(Y)$ can easily brought into upper triangluar form. 
By this method we find  				
\begin {equation}
det(M(Y)) = \varsigma_1 \cdot \prod_{b=2}^{r-1} \left(
\sum_{c=1}^b \varsigma_c + \sum_{c=a+1}^r \varsigma_c \prod_{e=b}^{c-1}q_e(z)^2
\right)
\cdot \prod_{a=2}^r \bar \varsigma_a
\end {equation}
with $q_e(z)$ as defined in equation \ref {E:qkz}.  Defining the 
$ r \times r $
matrix $\Delta $   as 
\begin {multline}
\Delta_{i,j} = 
(\delta_{i,j} -\delta_{i,j+1}) \cdot \left( \frac {1} {1-q_{i-1}(z)^2} \right) 
+ (\delta_{i,j} - \delta_{i,j-1}) \cdot \left( \frac {q_i(z)^2} {1- q_i(z)^2}
\right) \\
i,j = 1,\ldots,r; \quad q_0(z) \equiv q_r(z) \equiv 0
\end {multline}
and introducing the substitutions
\begin {equation}
\begin {split}
\bar \xi_i &:= \left( \frac {1} {1-x_i} \right) \cdot \bar \varsigma_i \\
\xi_i &:= \sum_{j=1}^r (\Delta^{-1})_{i,j} \varsigma_j
\end {split}
\end {equation}
with a little algebra, notably using 
\begin {equation}
\tilde U(Y) + \mathbf{I} = \Delta^{-1}
\end {equation}
($\mathbf{I}$ the unit matrix) one finds:
\begin {multline}
\sum_{w \in W_{2N}} P(N_{2k_1}(w),\ldots,N_{2k_r}(w)) \cdot z^{length(w)} = \\
z \frac {\partial}{\partial z} \Biggl[
\sum_{\substack{\vartheta \in D(k_1,\ldots,k_r) \\
g_i \in \N - \{0\}}}
\prod_{i=1}^r \hat K(x_i,k_{\vartheta(i)},h(0,1,z)) \\
\frac {(-1)} {(2\pi)^r} \int_{\R^{2r}} 
\left( \prod_{j=1}^r \frac {d\bar \xi_j\cdot d\xi_j} {1-x_j} \right)
\left( \prod_{k=2}^{r-1} \frac {\bar \xi_k \cdot \xi_k} {1 - q_k(z)^2} \right)
\frac {1} {1-q_1(z)^2} \\
\cdot \bar \xi_r \cdot 
\left( \xi_1 \left( \frac {1} {1-q_1(z)^2} \right) - \xi_2 \frac {q_1(z)^2}
{1 - q_1(z)^2} \right) \cdot exp \left( -\sum_{\substack{a=1 \\
b=1}}^r \bar \xi_a I_{a,b} \xi_b 
\right)
\Biggr]
\end {multline} 
with the matrix
\begin{multline}
I_{a,b} := \delta_{a,b} \cdot \left( 1 + x_a \cdot \left(
\frac {q_a(z)^2} {1 - q_a(z)^2} + \frac {q_{a-1}(z)^2} {1 - q_{a-1}(z)^2} 
\right) \right) \\
+ \delta_{a,b+1} \cdot \left( -x_{a-1} \frac {1} {1-q_{a-1}(z)^2} \right)
+ \delta_{a,b-1} \cdot \left( -x_{a+1}  \frac {q_a(z)^2} {1 - q_a(z)^2} \right)
\end{multline}
Expanding the exponential in the variables $x_a$ leaves us with a sum of Gaussian integrals over polynomials which can easily be performed. Applying the operators 
$\hat K$ we find 
\begin {multline}\label{E:mo1}
\sum_{w \in W_{2N}} P(N_{2k_1}(w),\ldots,N_{2k_r}(w)) \cdot z^{length(w)} = \\
z \cdot \frac {\partial}{\partial z} 
\begin{cases} T_2(z)^{(k_{1}),(k_{2})} \cdot (2 - \delta_{k_1,k_2}) & \Longleftarrow r =2 \\
\sum_{\vartheta \in D(k_1,\ldots,k_r)} & \\
\left< \Psi(k_{\vartheta(1)},z) \Bigr\rvert \prod_{j=2}^{r-2} Q(k_{\vartheta(j)},z)
\star \Phi(k_{\vartheta(r-1)},k_{\vartheta(r)},z) \right>_{z \diamond 2N} 
& \Longleftarrow  r > 2
\end {cases}
\end {multline}
with the notations of theorem \ref{T:one}  But from this theorem \ref{T:one} 
immediately follows, since a finite distribution is totally determined 
by it's moments. 
\begin {remark}
From equation \ref{E:psid} and \ref{E:qdef}
the matrix elements $Q_{(\varrho,t),(\tilde \varrho,\tilde t)}(j,z)$ and the
components $\Psi_{(\varrho,t)}(j,k,z)$ are $0$ for 
$\varrho + t \geq j - 1$. Therefore for $\varrho + t \geq j - 1$  we have
$(Q(j,z)^n)_{(\varrho,t),(\tilde \varrho,\tilde t)} = 0 \quad \forall n \geq 1$. Therefore 
for calculating the distribution of the $k$-multiple point range whose multiplicity is bound by
$2\cdot k_{max}$ with $k_{max} \in \N $ the matrices/vectors $Q$, $\Psi$ and therefore also $\Phi$ only
need to be evaluated in the $(k_{max} - 1)^2$ dimensional subspace in which the 
elements of $Q$ and $\Psi$ are different from zero. This is independant of $N$. 
\end {remark}
\begin{example}
From the above remark we can immediately calculate the following examples for 
singlepoints ($A(z)$ and $B(z)$ see theorem \ref{T:one})
\begin{eqnarray}\label{E:dsin}
\sum_{\substack{w \in W_{2N} \\ N_2(w) = 0}} z^{length(w)} 
& = & 
\left( A(z) - 1 
 + z \frac {\partial}{\partial z} \left( A(z)^2 \sum_{f=1}^{\infty} f \cdot \frac {B(z)^{2f}} {1 - B(z)^{2f}} \right)
\right)_{z \diamond 2N} \nonumber \\
\sum_{\substack{w \in W_{2N} \\ N_2(w) = 1}} z^{length(w)} 
& = & 
\left( \frac{4z^2}{A(z)}
 - 2 \cdot z \frac {\partial}{\partial z} 
\left( A(z)^2 \sum_{f=1}^{\infty} f \cdot \frac {B(z)^{2f}} {1 - B(z)^{2f}} \right)
\right)_{z \diamond 2N} \\ 
\sum_{\substack{w \in W_{2N} \\ N_2(w) = 2}} z^{length(w)} 
& = & 
\left( 
 z \frac {\partial}{\partial z} 
\left( A(z)^2 \sum_{f=1}^{\infty} f \cdot \frac {B(z)^{2f}} {1 - B(z)^{2f}} \right)
\right)_{z \diamond 2N} \nonumber
\end{eqnarray}
and doublepoints
\begin{multline}\label{E:ddou}
\sum_{w \in W_{2N}} e^{it\cdot N_4(w)}\cdot z^{length(w)} = \Biggl(
\frac {1}{A(z)} - 1 + (e^{it} - 1)\cdot \frac {4z^2 \cdot (1 - A(z))} {A(z)} 
+  \\ 
z\cdot \frac{\partial}{\partial z} \Biggl(
A(z)^4 \cdot (e^{it} - 1)^2
 \sum_{f=1}^{\infty} \frac{1}{f} \left( \frac {f \cdot (1 - A(z))}
{A(z)} - \binom {f}{2}
\right)^2 \frac {B(z)^{2f}} {1 - B(z)^{2f}} \\ 
+ A(z)^6 \cdot (e^{it} -1)^3 
\cdot \left( 
\sum_{f = 1}^{\infty} \left( 
\frac {f}{A(z)} - \frac {f \cdot (f + 1)}{2}
\right)
\frac {B(z)^{2f}}{1 - B(z)^{2f}}
\right)^2
\cdot  \\
\frac {1} {1 - (e^{it} - 1)\cdot A(z)^2 \cdot  \sum_{f=1}^{\infty}
f \cdot \frac{B(z)^{2f}} {1 - B(z)^{2f}}}
\Biggr) \Biggr)_{z \diamond 2N}
\end{multline}
\end{example}
\begin {example}
As an example for the moments (equation \ref{E:mo1}):
\begin{multline}
\sum_{w \in W_{2N}} \binom {N_2(w)}{2} \binom {N_6(w)}{2} z^{length(w)}
=  \\ 2 z \frac {\partial}{\partial z} \Biggl[
2 \left( G_{1,1}(z)^2 (1 - A(z)) - 2G_{1,2}(z)G_{1,1}(z) \right) \cdot  \\
\left( G_{1,1}(z) (1 - A(z)) - G_{1,2}(z)) \right) + 
G_{1,1}(z)^2H_{2,4}(z) -  \\  
2 G_{1,1}(z) \left(
G_{1,2}(z) G_{1,1}(z) (1 - A(z)) - G_{1,1}(z)G_{1,3}(z) - G_{1,2}(z)^2
\right)
\Biggr]_{z \diamond 2N}
\end{multline}
\end {example}
\begin {example}
The distribution of the range of a closed onedimensional random walk can be 
obtained by summing over the multiplicities in the equations of theorem 
\ref{T:one}. The summation greatly simplifies the formulas and after some algebra
one finds
\begin {multline}\label{E:dvol}
\sum_{\substack{w \in W_{2N} \\ ran(w) = m}}  z^{length(w)} = \\
z \cdot \frac {\partial}{\partial z}\Biggl(
\ln(1-B(z)^{2(m-1)}) - 2\ln(1-B(z)^{2m})  + \ln(1-B(z)^{2(m+1)})
\Biggr)_{z \diamond 2N}
\end {multline}
\end {example}
\section{ The asymptotic behaviour}\label{S:asy}

\begin {notation}
In this chapter we denote by the symbol $\sqrt{\quad}$ the function: 
\begin{eqnarray}
\sqrt{\quad}:& \C- \bigr]-\infty,0 \bigr] & \longmapsto  \C \\
     &\varsigma & \longmapsto  \sqrt{\varsigma} \nonumber
\end{eqnarray}
the branch of the square root which is nonegative on the part 
of the real axis where it is defined. 
\end {notation}
\begin {lemma}
For any $\varsigma \in \C - \bigl[1,\infty \bigl[ $   and 
$k \in \N - \{0\}$  
the infinite sum
\begin{equation}
\begin{split}
\sum_{f=1}^{\infty} f^k \frac {b(\varsigma)^f} {1 - b(\varsigma)^f} &=
\sum_{f=1}^{\infty} \sigma_k(f) \cdot b(\varsigma)^f \\
b(\varsigma) &:= \frac {1- \sqrt{1-\varsigma}}{1 + \sqrt{1 - \varsigma }}
\end{split}
\end {equation}
converges absolute. ( $\sigma_k(f)$ the sum of the $k$'th powers of the divisors 
of $ f $ \cite{poc})
\end {lemma}
\emph {Proof.}
From expressing $\sqrt{\quad}$ in polar coordinates it is immediately 
clear that the real part of $\sqrt{\quad}$ is strictly positive on 
$ \C - \bigr]-\infty, 0 \bigr]$.
Therefore 
\begin {equation}\label{E:abs}
\mid b(\varsigma) \mid = \frac
{\sqrt{(1-\Re (\sqrt{1-\varsigma}))^2 + \Im (\sqrt{1-\varsigma})^2}}
{\sqrt{(1+\Re (\sqrt{1-\varsigma}))^2 + \Im (\sqrt{1-\varsigma})^2}} < 1
\end {equation}
($\Re$ and $\Im$ denoting the real and imaginary 
part) for $\varsigma \in \C - \bigl[1, \infty \bigl[$ . 
As it is obvious that $\mid \sigma_k(f) \mid < f^{k+1} $  the lemma is proven. 
\begin{corollary}
For $k \in \N - \{0\}$  the function 
\begin{eqnarray}
g_k:& \C- \bigl[1,\infty \bigl[ & \longmapsto \C \nonumber \\
    & \varsigma & \longmapsto (\sqrt{1 - \varsigma})^{k+1} \sum_{f=1}^{\infty}
f^k \frac {b(\varsigma)^f}{1 - b(\varsigma)^f} 
\end{eqnarray}
is a holomorphic function and therefore for $\mid \varsigma \mid < 1$  there is 
a Taylor expansion
\begin {equation}
g_k(\varsigma) = \sum_{n=0}^{\infty} c_n(k) \cdot \varsigma^n
\end {equation}
\end {corollary}   
From the theory of power series it is well  known that the absolute 
value of the Taylor coefficients of a power series $\wp $ with radius 
of convergence $1$ are $O(1/n^m)$ if $\wp $ together with it's derivatives 
up to $m$'th order has a contiuous contination onto $\bar U_1(0) :=
\{ \varsigma \in \C \rvert \mid \varsigma \mid  \leq 1 \} $. As the 
functions $g_k$ are regular on $\bar U_1(0) - \{1\}$
the asymtotic behaviour of $c_n(k)$  for $n \longmapsto \infty $  can be 
obtained by a study of the singularity at $\varsigma = 1$. To obtain this 
behaviour we use the 
\begin{theorem}\label{T:asy}
For any $M \in \N - \{0,1\} $ there are (easily calculatable) constants $ \lambda_l(k),\tilde\lambda_l(k)$  
such that the functions 
\begin{eqnarray} 
\nonumber 
\wp_k : & U_1(0) & \longmapsto  \C \\
& \varsigma  & \longmapsto   g_k(\varsigma) - \frac {k! \cdot \zeta(k+1)} {2^{k+1}} + \sum_{l = 1}^M \lambda_l(k) \cdot (1 - \varsigma)^l \nonumber \\
& & \qquad + \delta_{k,1} \cdot \sqrt{1-\varsigma} \cdot \left( \frac{1}{4} + \sum_{l=1}^M \tilde \lambda_l(k)\cdot (1 - \varsigma)^l \right)
\end{eqnarray}
together with their derivatives in $ \varsigma $  up to the order $M$ have 
a continuous continuation onto  $ \bar U_1(0)$ which is 0 for $\varsigma = 1$ in all these orders.
\end {theorem}
\emph{Proof.}
Because of equation \ref{E:abs} we can concentrate on a pierced 
vicinity $V$ of $\varsigma = 1$  in $\bar U_1(0) - \{ 1 \}$ . We choose  
$m > M$
and interprete the sums
\begin{equation}
\frac{g_k(\varsigma)}{\left( \sqrt{1- \varsigma} \right)^{k+1}} = 
\sum_{f = 0}^{\infty} f^k \frac {b(\varsigma)^f} {1 - b(\varsigma)^f}
\end{equation}
as trapezoidal sums of the integrals 
\begin{equation}
\int_0^{\infty} \frac {y^k} {b(\varsigma)^{-y} - 1} dy
\end{equation} 	
According to the Euler-Mac Laurin expansion of the trapezoidal 
sum \cite{stoer} for step width $1$ we can write:
\begin {multline}\label{E:eul}
\frac {g_k(\varsigma) } {\left( \sqrt{1 - \varsigma} \right)^{k+1} }
= \sum_{f = 0}^{\infty} f^k \frac {b(\varsigma)^f} { 1 - b(\varsigma)^f} =
\\
\frac {\delta_{k,1}} {2 \cdot \ln(b(\varsigma))} 
+ \frac {(-1)^{k+1}} {\ln(b(\varsigma))^{k+1}} \left(
\zeta(k+1) \cdot k! - \sum_{j \geq 1 \wedge 2j \geq k}^m \frac{B_{2j} \cdot B_{2j-k} \cdot \left( \ln(b(\varsigma))\right)^{2j}}{2j*\left(2j-k \right)!} \right)\\
- \frac {1}{(2m+2)!} \int_0^{\infty} \left( K_{m} (y) - K_{m}(0) \right) \cdot
 \left( \left( \frac {\partial }{\partial x} \right)^{2m+2}
\left( \frac {x^k} {b(\varsigma)^{-x} - 1} \right) \Biggr\rvert_{x = y} \right)dy  
\end{multline}
with the function 
\begin{eqnarray}
K_{m}: & \R & \longmapsto \R \\
&x &\longmapsto B_{2m+2}(x-i) \Longleftarrow x \in \left[ i,i+1 \right] \nonumber
\end{eqnarray}
and $ B_{j}(x) $ the $j$'th Bernoulli Polynomial and the Bernoulli numbers $B_j = B_j(0)$. The function $K_{m}$ is periodic 
with period $1$ and bound on $\R$. From equation \ref{E:eul} it is clear that 
theorem \ref{T:asy} is true (and how the constants $\lambda_l,\tilde \lambda_l$ are obtained), if  the integral term on the right hand 
side multiplied with $\left( \sqrt{1-\varsigma} \right)^{k+1}$  
(together with it's appropriate derivatives 
in $\varsigma $ ) has a continuous continuation onto $\bar U_1(0) $ which is 0 for 
$\varsigma = 1$. We split the integral 
in two parts, one $( I_1 )$ over the interval $\left[0,1/\vert\ln(b(\varsigma))\vert \right]$  
the other one $( I_2 )$ over the interval $\bigl[1/\vert \ln(b(\varsigma))\vert,\infty \bigr[$. 
For $I_1$ we use the  well known Taylor expansion \cite{dan358}
\begin{multline}
\varphi_{m,k}(x,\varsigma) := \left( \frac{\partial} {\partial x} \right)^{2m+2} 
\left( \frac {x^k} {b(\varsigma)^{-x} - 1 }\right) = 
\left( \ln(b(\varsigma)) \right)^{2m+2 - k} \cdot \\
\sum_{n = \max(0,k - 2m -3)}^{\infty} \frac{B_{n +2m + 3 - k}} {(n+2m +3 -k)!}
\binom {n + 2m + 2} {2m + 2} \left(-\ln(b(\varsigma)) \right)^n \cdot x^n
\end{multline}
The inequality 
\begin{equation} \label{E:ine}
\Biggl\lvert \left( \frac {\partial} {\partial \varsigma} \right)^s
\left( \ln(b(\varsigma)) \right) \Biggr\rvert <
C_s \cdot 2^n \cdot n^s \lvert \left( \sqrt{1 - \varsigma} \right)^{-2s} 
\ln(b(\varsigma))^n \rvert
\end{equation}
( $C_s$ real constants) holds independantly of $n$ and $\varsigma$ in a pierced 
vicinity of $\varsigma = 1$   in $\C$ because of the Taylor expansion of the logarithm 
and $b(\varsigma)$. But from this it can easily be seen that 
\begin{equation}\label{E:i1}
\Biggl\lvert \left( \frac {\partial}{\partial \varsigma} \right)^s I_1 \Biggr\rvert <
C_{s,m} \Bigl\lvert \sqrt{1- \varsigma} \Bigr\rvert^{2m + 2 - 2s -k -1}
\end{equation}
with certain real constants  $C_{s,m}$ in a pierced vicinity of $\varsigma = 1$ in 
$ \C $. 
For $I_2$ we realize that for $\mid \varsigma \mid \leq 1 $
\begin{equation}
\Re \left( \sqrt{1-\varsigma} \right) \geq \mid \Im \left( \sqrt{1 - \varsigma} \right) \mid
\end{equation}
and therefore 
\begin{equation}\label{E:i2}
\begin{split}
\Re \left( -\ln(b(\varsigma)) \right)& > \frac{1}{\sqrt{5}} \cdot \mid \ln(b(\varsigma)) \mid \\
\Biggl\lvert \frac {1}{b(\varsigma)^{-y} - 1} \Biggr\rvert & < 
\frac {1}{1 - e^{-\frac{1}{\sqrt{5}}}} \cdot 
e^{- \frac{y \mid \ln(b(\varsigma)) \mid}{\sqrt{5}}} \quad 
\text{for $y > \frac{1}{\mid \ln(b(\varsigma))\mid} $}
\end{split}
\end{equation}
in a pierced vicinity of $\varsigma = 1$  in $\bar U_1(0)$. 
With elementary calculations it can be shown that
\begin{displaymath}
\left( \frac {\partial}{\partial \varsigma} \right)^s \varphi_{m,k}(x,\varsigma)
\end{displaymath}
can be written as  $1/( b(\varsigma)^{-x} - 1 )$ times a polynomial in the 
variables
\begin{equation}
\left( \frac {\partial} {\partial \varsigma} \right)^j (-\ln(b(\varsigma )), \;
\frac{1}{\left( b(\varsigma)^{-x} - 1\right)}, \;
\frac{b(\varsigma)^{-x}}{\left( b(\varsigma)^{-x} - 1\right)}, \;
x
\end{equation}
$j = 0,\ldots,s$ with the features 
\begin{enumerate}
\item that if we add up the sum of the derivatives before the factors 
$\left( -\ln(b(\varsigma)) \right)$  we always get $ s $
\item that the degree of each term in x differs from the number 
of all factors $\left( -\ln(b(\varsigma)) \right)$ (with or without a 
derivative in front of them)  by $2m + 2 - k$. 
\end{enumerate}
But with equation \ref{E:i2} that means that 
\begin{equation}
\Biggl\lvert \left( \frac {\partial} {\partial \varsigma} \right) I_2 \Biggr\rvert
< \tilde C_{s,m} \lvert \sqrt{1- \varsigma}\rvert^{2m+2 -2s -k -1}
\end{equation}
 with certain real constants $\tilde C_{s,m}$   in a pierced vicinity of $\varsigma = 1$ in 
$\C$ which together with equation \ref{E:i1} proves the theorem.
\begin{example}
In this and the following examples let $Pr_n()$ denote the probability for closed simple random walks of length $2n$ and $E_n()$ the corresponding expectation values as in 
equation \ref{E:amo}. 
Then 
\begin{equation}\label {E:dou}
\begin{array}{rcl}
Pr_n(N_2(w) = 0 ) & = &\frac{1}{4}  - \frac{1}{4n} - \frac{1}{48n^2} + \frac{13}{144n^3} + \frac{421}{2880n^4} + O(\frac{1}{n^5}) \\ \\
Pr_n(N_2(w) = 1 ) & = &\frac{1}{2} - \frac{5}{24n^2} - \frac{11}{36n^3} - \frac{511}{1440n^4}   + O(\frac{1}{n^5}) \\ \\
Pr_n(N_2(w) = 2 ) & = &\frac{1}{4} + \frac{1}{4n} + \frac{11}{48n^2} + \frac{31}{144n^3} + \frac{601}{2880n^4} + O(\frac{1}{n^5}) \\ \\
Pr_n(N_4(w) = l > 2 ) & = &\alpha^l \cdot \left( \theta_0 + l \cdot \theta_1  \right) + O(\frac{1}{n})
\end{array}
\end{equation}
with 
\begin{equation}
\begin{split}
\alpha & := \frac{\pi^2}{24 + \pi^2} \\
\theta_0 & := \frac{216}{\pi^6} \frac {\left(\frac{\pi^2}{3} - \zeta(3)\right)^2}
{1 + \frac{\pi^2}{24}} \cdot 
\left( 
\frac{4 + \frac{\pi^2}{3}    }{\frac{\pi^2}{3} - \zeta(3)} - 
\frac{3 + \frac{\pi^2}{24}    }{\frac{\pi^2}{6}\cdot \left( 1 + \frac{\pi^2}{24} \right)} -
\frac{3}{4 \cdot \left( 1 + \frac{\pi^2}{24}    \right)}
\right)\\
\theta_1 & := \frac {1296}{\pi^8} \cdot \left( 
\frac {\frac{\pi^2}{3}  - \zeta(3)  }{1 + \frac{\pi^2}{24}}
\right)^2
\end{split}
\end{equation}
\end{example}
\begin{example}
From the following table about the number of doublepoints of a onedimensional 
closed simple random walk of length $78$ one gets a feeling both of the distribution 
$Pr_n(N_4(w) = l)$ for $n = 39$
and of the accuracy of the asymptotics $\Pr_{\infty}(N_4(w) = l) := \lim_{n \rightarrow \infty} Pr_n(N_4(w) = l)$
as given by equation \ref{E:dou} for $l > 2$. For $l \leq 2$ the results were 
calculated from the corresponding formulas not give here to save space. 
\begin{center}
\textbf{Table 1: Doublepoints}
\end{center}
\begin{tabular}{c|r|c|c}
\hline
$l$ & 
$\sum_{\substack{ w: length(w) = 78 \\ N_4(w) = l}} 1$ & 
$Pr_{39}(N_4(w) = l)\approx $ &
$Pr_{\infty} (N_4(w) = l) \approx $ \\
\hline
0 & 9379489746558670340000 & 0.34462 & 0.35101 \\
\hline
1 &11080781119308072700000 & 0.40713 & 0.40526 \\
\hline
2 & 4768982388008920550000 & 0.17522 & 0.17199 \\
\hline
3 & 1321976178995539300000 & 0.04857 & 0.04779 \\
\hline
4 &  446940016375442637000 & 0.01642 & 0.01608 \\
\hline
5 &  148016854282117480000 & 0.00544 & 0.00531 \\
\hline
10 &    478890500239691072 & 0.0000176 & 0.0000177 \\
\hline
\end{tabular}
\end{example}
\begin{example}
For the distribution of the multiple range with higher multiplicity $k>2$ we find the general form
\begin{multline}
Pr_n(N_{2k}(w) = l > 2 )  = \sum_{i=1}^{k-1}  \alpha_i(k)^l \cdot \left( \theta_0(k)^{[i]} + l \cdot \theta_1(k)^{[i]} \right)  \\+ \sum_{1\leq i < j \leq k-1}\theta_2(k)^{[i][j]} \cdot \sum_{m=1}^{l-1}\alpha_i(k)^m \cdot \alpha_j(k)^{l-m}  +  O(1/n)
\end{multline}
with certain complex constants $\alpha_i(k),\theta_0(k)^{[i]}, \theta_1(k)^{[i]}, \theta_2(k)^{[i][j]}$. To give a feeling for the distribution for large $l$, below we present a list of the numerical values of the constants $\alpha_i(k)$ with $ 1\leq i \leq 6$ sorted according to their absolute value. 
\begin{center}
\textbf{Table 2: Dominating constants for higher multiplicity}
\end{center}
\begin{center}
\begin{tabular}{c|c|c|c|c|c|c}
\hline
$k$ & 
$\alpha_1(k)\approx $ & $\alpha_2(k)\approx $ & $\alpha_3(k)\approx $ & $\alpha_4(k)\approx $ & $\alpha_5(k)\approx $ & $\alpha_6(k)\approx $ \\
\hline
2 & 0.29140 & & & & & \\
\hline
3 & 0.29018 & -0.23057 & & & & \\
\hline
4 & 0.29867 & -0.14176 & 0.12556& & & \\
\hline
5 & 0.30263 & -0.18822 & -0.08169 & 0.07648 & & \\
\hline
6 & 0.30590 & -0.13810 & 0.11967 & -0.04539 & 0.04382 & \\
\hline
7 & 0.30829 & -0.16093 & -0.09351 & 0.08564 & -0.02471& 0.02426\\
\hline 
8 & 0.31020 & -0.12830 & 0.11074 & -0.05984 & 0.05671 & -0.01327\\
\end{tabular}
\end{center}
\end{example}
\begin{example}
To give an impression of the mutual dependancy among the variables $N_{2k}(w)$ we also give asymptotic formula for the second moments:
\begin{multline}\label {E:2mom}
E_n(N_{2k_1}(w) \cdot N_{2k_2}(w)) = \\
\delta_{k_1,k_2} + \frac {1}{2} + 2 \cdot \sum_{r_1 = 1}^{k_1} \sum_{r_2 = 1}^{k_2}
\binom {k_1 - 1}{r_1 - 1} \cdot \binom {k_2 - 1} {r_2 -1} \cdot (-1)^{r_1 + r_2}
\cdot \binom {r_1 + r_2} {r_1} \frac {1}{2^{r_1 + r_2}} \cdot \\
\left( 
\frac {k_1 + k_2 - r_1 - r_2} {r_1 + r_2} \zeta(r_1 + r_2) + 
\frac{\binom{r_1}{2} + \binom{r_2}{2}}{\binom{r_1 + r_2}{2}} \cdot 
\begin{cases} 0 & if \> r1 + r2 = 2 \cr \zeta(r_1 + r_2 - 1) &else 
\end{cases}
\right) \\ + O(1/n)
\end{multline}
In numbers this formula means
\begin{center}
\textbf{Table 3: $\lim_{n \rightarrow \infty}(E_n(N_{2k}(w)\cdot N_{2\tilde k}(w))-E_n(N_{2k}(w))\cdot E_n(N_{2\tilde k}(w))) \approx $}
\end{center}
\begin{center}
\begin{tabular}{c|c|c|c|c|c|c}
\hline
$k/\tilde k$ & 1 & 2 & 3 &4 &5 & 100 \\
\hline
1& 0.50000 & -0.08877 & 0.02195 & 0.03509 & 0.02398 & 0.00000 \\
\hline
2&         & 1.02195  & 0.10274 & 0.12970 & 0.12418 & 0.00500 \\
\hline
3&         &          & 1.16494 & 0.19628 & 0.19974 & 0.01000 \\
\hline
4&         &          &         & 1.23626 & 0.251235 & 0.01500 \\
\hline 
5&    		 &          &         &         & 1.27949  & 0.02000 \\
\hline 
100 &      &          &         &         &          & 1.47074 \\
\hline 
101 &      &          &         &         &          & 0.47061 \\
\end{tabular}
\end{center}
\end{example}
\begin{example}
For the moments of the range we get the asymptotic formula
\begin{equation}\label{E:avol}
E_n(ran(w)^r) =
= \xi(r) \cdot (E_n(ran(w)))^r \cdot (1 + O(1/n))
\end{equation}
with the $\xi$ function of Riemann 
\begin{equation}
\xi(r) = r \cdot (r-1) \cdot \zeta(r) \cdot \Gamma(\frac{r}{2})\cdot \pi^{- \frac {r}{2}}
\end{equation}
\end{example}

\section{Summary}\label {S:out}
In this article the joint distribution of the $k$-multiple point range in a closed simple random walk was  discussed. It was shown how to obtain the moments of the distribution by studying the 
number of functions between Eulerian graphs and walks preserving basic structural features. From this study 
\begin {itemize}
\item the first moment of the distribution of points of any multiplicity in any 
dimension
\item the exact joint generating function in dimension one
\item a general method to calculate the asymptotic behaviour of moments and finite subdistributions for walks of large lengths in dimension one
\end {itemize}
were derived. The joint distribution in one dimension turns out ot be a generalized
geometric distribution. 

The results in the paper can probably be easily extended to random walks of mean
zero and finite variance and without the restriction of being closed without much effort.
Mainly the functions $h(x,d,z)$ will change in the corresponding equations. For the 
onedimensional case the mutliplicative structure of $h(x,1,z)$ as expressed in equation
\ref {E:qkz} is the crucial point for the proof, so ondimensional random walks with
this property can at once be solved in the same way. 

Obviously the distribution in higher dimensions is of even greater interest than in one dimension. It will be discussed in articles to follow. The limit theorems obtained in the literature \cite{Ham1, Ham2, Ham3} give an 
asymptotic picture of these distributions. Calculating them for an ensemble
of walks instead of one walk would open the possibility for a more rigorous definition of and calculation in the path integral formulation of an important class of quantum field theories \cite{bry}. The interesting field theories are nontrivial,  
the central limit theorems in the language of field theory should give the 
mean field results.  
A more refined analysis of higher moments is necessary for a better understanding
of QFT. 

The moments of the asymptotic distribution of the range of a 
onedimensional walk are related to Riemanns $\xi$-function. So the Riemann conjecture can
be reinterpreted as a statement about the zero points of the characteristic function of the rescaled asymptotic distribution of the range of the onedimensional random walk.


\begin{thebibliography}{9}
\bibitem{Mard} 
Mardudin, A.A., Montrose, E.W., Weiss, G.H., Hermann, R., Milne, A.
Acad. Roy. Belg. Cl. Sci. Mem. 4, 1960, pp. \ 147ff
\bibitem{Tek} Dvoretzky, A., Erd\"os P., Kakutani S., 
Acta Sci. Math. Szeged 12, 1950, pp. \ 75ff.
\bibitem{Pitt} Pitt, J.H., ``Multiple points of transient random walks,'' Proc. Amer. Math. Soc. 43, 1974, pp. \ 195--199
\bibitem{Ham1} Hamana, Y., ``On the central limit theorem for the multiple point range
of random walks,'' J. Fac. Sci. Univ. Tokyo 39, 1992, pp. \ 339--363
\bibitem{Ham2} Hamana, Y., ``On the multiple point range of three dimensional
random walks,'' Kobe J. Math., 12, 1995, pp. \ 5--122
\bibitem{Fla} Flatto, L., ``The multiple range of the two dimensional recurrent walk,''
Ann. Probab. 4, 1976, pp. \ 229--248
\bibitem{Ham3} Hamana, Y., ``A remark on the multiple point range of two-dimensional
random walks,'' Kyushu J. Math. 52, 1998, pp. \ 23--80
\bibitem{Jain} Jain, N.C., Pruitt, W.E., ``The range of random walks,'' Proc. Sixth 
Berkeley Symp. Math. Stat. Probab., Berkeley, CA, 1973, pp. \ 31--50
\bibitem{QFT} Zinn-Justin, J., ``Quantum Field Theory and Critical Phenomena,''
Clarendon Press, Oxford 1989, pp. \ 4--5
\bibitem{aar}van Aardenne-Ehrenfest, T., de Brujin, N.G., Simon Stevin 
 Wis. Natuurkd. Tijdschr. 28, 1951, pp. \ 203ff.
\bibitem{poc} Danos, M., Rafelski, J., ``Pocketbook of Mathematical Functions,''
Harri Deutsch, Thun, Frankfurt/Main, 1984, pp. \ 370
\bibitem{stoer} Stoer, J., ``Einf\"uhrung in die Numerische Mathematik,'' Springer,
New York, Heidelberg, Berlin, 1979, vol. 1, pp. \ 104--109
\bibitem{dan358} Danos, M., Rafelski, J., ``Pocketbook of Mathematical Functions,''
Harri Deutsch, Thun, Frankfurt/Main, 1984, pp. \ 358
\bibitem{bry} Brydges, D., Fr\"ohlich, J., Spencer, T., ``The Random Walk Representation 
of Classical Spin Systems and Correlation Inequalities,''  Commun. Math. Phys., 83,
1982, pp. \ 123--150
\end{thebibliography}
\end{document}